\documentstyle[12pt]{article}

\input{amssym.def}
\input{amssym.tex}
\begin{document}
\title{Geodesic laminations revisited}

\author{Igor  ~Nikolaev\\
Department of Mathematics\\
2500 University Drive N.W.\\
Calgary T2N 1N4 Canada\\
{\sf nikolaev@math.ucalgary.ca}}


 \maketitle


\newtheorem{thm}{Theorem}
\newtheorem{lem}{Lemma}
\newtheorem{dfn}{Definition}
\newtheorem{rmk}{Remark}
\newtheorem{cor}{Corollary}
\newtheorem{prp}{Proposition}
\newtheorem{exm}{Example}


\newcommand{\N}{{\Bbb N}}
\newcommand{\F}{{\cal F}}
\newcommand{\R}{{\Bbb R}}
\newcommand{\Z}{{\Bbb Z}}
\newcommand{\C}{{\Bbb C}}

\begin{abstract}
The Bratteli diagram is an infinite graph which reflects the
structure of projections in a $C^*$-algebra. We prove that
every strictly ergodic unimodular Bratteli diagram of rank $2g+m-1$
gives rise to a minimal geodesic lamination with the $m$-component  principal region
on a surface of genus $g\ge 1$. The proof is based on the Morse theory
of the recurrent geodesics on the hyperbolic surfaces.

\vspace{7mm}

{\it Key words and phrases: Bratteli diagrams, 
geodesic laminations}

\vspace{5mm}
{\it AMS (MOS) Subj. Class.:  19K, 46L, 57M.}
\end{abstract}

\section*{Introduction}
The paper deals with three apparently  independent topics: the Morse theory
of the recurrent geodesics, the Nielsen-Thurston theory of the geodesic 
laminations and, finally, a piece of the $C^*$-algebra theory, known as
the Bratteli diagrams. Our goal is to show that 
the Bratteli diagrams imply   the geodesic laminations   
via the Morse theory. Such a result links  geometry to the operator
algebras.

Recall that
the simplest $C^*$-algebra can be written as
\begin{equation}
A_n=M_{n_1}({\Bbb C})\oplus\dots\oplus M_{n_r}({\Bbb C}),
\end{equation}
where $M_{n_i}$ are the square matrices with the complex entries and $r$, $n_1,\dots, n_r$ some 
non-negative integers. The closure of an infinite sequence
\begin{equation}
A_1\subseteq A_2\subseteq\dots 
\end{equation}
of such $C^*$-algebras is a $C^*$-algebra denoted  by
\begin{equation}
{\goth A}=\lim_{n\to\infty}A_n.
\end{equation}

The Bratteli diagram is an infinite connected graph, which reflects an  embedding of the 
finite-dimensional algebras $A_n$ in the inductive limit ${\goth A}$ (\cite{Bra}).  
Such an  embedding is described by a matrix of partial multiplicities $M_n$ whose entries
are non-negative integers. If $Rank ~M_n=r=Const$ through all the embeddings, we
say that the Bratteli diagram has rank $r$. If $det~M_n=\pm 1$ for all $n$, the Bratteli
diagram is called {\it unimodular}.  The Bratteli diagram is {\it strictly ergodic} 
if the linear space
\begin{equation}
L=\bigcap_{n=1}^{\infty} M_1\dots M_n ~({\Bbb R}^r_+)
\end{equation}
has dimension $1$.

Let $S$ be a connected complete hyperbolic surface of genus $g$. Recall that a {\it geodesic} on $S$
is the maximal arc consisting of locally shortest sub-arcs. 
The study of geodesics on surfaces goes back to Birkhoff, Hadamard, Morse and Hedlund.  
A geodesic is {\it simple} if it has no self-crossing or self-tangent points.
The periodic geodesic is an elementary example of simple geodesic. 
Birkhoff conjectured and Morse proved existence of simple recurrent non-periodic
geodesics. The set of such geodesics turns to be uncountable on any hyperbolic surface.

The problem of classification of the recurrent non-periodic geodesics
leads to the concept of a {\it geodesic lamination}. 
The geodesic lamination $\lambda$ on $S$ is a disjoint union  of all recurrent
non-periodic geodesics which lie in the closure of each other.
One can think of $\lambda$ as an uncountable  set of non-periodic geodesics
running in the ``same direction'' on surface $S$.
The geodesic lamination proved to be fundamental in the dynamics, complex analysis,
and topology (\cite{CaB}).

The set $S-\lambda$ is a {\it principal region} of the geodesic lamination $\lambda$. 
It is an important combinatorial invariant of $\lambda$ and essentially a finite union
of the ``ideal polygons'' $U_1,\dots,U_m$. Each $U_i$  has a {\it type}
described by a positive integer $k_i$. 
By a {\it singularity data} of $\lambda$
one understands a set $\Delta=(k_1,\dots,k_m)$ such that $\Sigma k_i=2g-2$. 
In the present paper we prove the following theorem.
\begin{thm}\label{thm1}
Let $B$ be a strictly ergodic unimodular  Bratteli diagram of rank $r\ge 2$
and $\Delta$ a singularity data. Then the pair $(B,\Delta)$ 
 defines  a  geodesic lamination $\lambda$ on the hyperbolic
surface $S$ of the genus $g=g(\Delta)$ such that:

\medskip
(i) the principal region of $\lambda$ has $r-2g+1$ connected components;

\smallskip
(ii) the singularity data of $\lambda$ coincides with $\Delta$. 
\end{thm}

\medskip\noindent
The paper is organized as follows. In Section 1 we introduce the
notation and lemmas which will be used to prove our main theorem.
For other facts, we refer to the relevant bibliography.  
Theorem \ref{thm1} is proved in Section 2. An example of 
the Bratteli diagram of an irrational rotation algebra with
the ``golden mean''  Rieffel parameter is considered in Section 3.

\medskip\noindent
{\bf Acknowledgments} It is my pleasure to thank F.~Bonahon, O.~Bratteli,  
and  G.~A.~Elliott for  helpful 
discussions of the subject of this note.

\section{Preliminaries}
In this section we bring together some useful facts on the
geodesic laminations, Morse theory and combinatorics of 
the $C^*$-algebras. Our exposition is sketchy, and we refer 
the reader to (\cite{Bra}), (\cite{CaB} and (\cite{Mor}) for
a complete treatment. We omit an introduction to  the interval exchange
transformations, which are technically important for the proof of Theorem \ref{thm1}. 
 However, we hope that the reader can recover the details by reading (\cite{Vee1})
and (\cite{Vee2}). Our notation is borrowed from there.

\subsection{Geodesic laminations}
Let $S$ be a connected complete hyperbolic surface. By a {\it geodesic}
on $S$ we understand the maximal arc consisting of locally shortest sub-arcs.
The geodesic is called {\it simple} if it has no self-crossing or self-tangent
points.

A {\it geodesic lamination} on $S$ is a closed subset $\lambda$ of $S$
which is a disjoint union of simple geodesics of $S$. The geodesics
of $\lambda$ are called {\it leaves} of $\lambda$. The lamination
$\lambda$ is called {\it minimal} if no proper subset of $\lambda$
is a geodesic lamination. The following lemma gives a classification of
the minimal laminations.
\begin{lem}\label{prp1}
A minimal lamination $\lambda$ in a closed orientable hyperbolic
surface $S$ is either a singleton (simple closed geodesic) or an
uncountable nowhere dense subset of $S$.
\end{lem}
{\it Proof.} See Lemma 3.3 of
(\cite{CaB}).
$\square$

\bigskip
The laminations $\lambda,\lambda'$ on $S$ are (topologically)
{\it equivalent} if there exists a homeomorphism $\varphi:S\to S$
such that each leaf of $\lambda$ through point $x\in S$ goes
to the leaf of $\lambda'$ through point $\varphi(x)$. Clearly,
the set of all laminations, $\Lambda(S)$, on $S$ splits
into the equivalence classes under this relation.

\bigskip
If $\lambda\in\Lambda(S)$, then a component of $S-\lambda$ is called
a {\it principal (complementary) region} for $\lambda$.
(Note that $S-\lambda$ may have several connected components.)
The leaves of $\lambda$ which form the boundary of a principal
region are called {\it boundary leaves.}
If $\lambda$ is minimal (which we always assume to be), then each boundary
leaf is a dense leaf of $\lambda$, isolated from one side.

Note that by Lemma  \ref{prp1}
\begin{equation}\label{eq3}
Area~(S-\lambda)=Area~S,
\end{equation}
and therefore principal region is a complete hyperbolic surface
of area $-2\pi\chi(S)$, where $\chi(S)=2-2g$ is the Euler characteristic
of surface $S$.  If $U$ is a component of the preimage of
the principal region in $D$, then $U$ is the union of
the ideal polygons $U_i$ in $D$ (see Fig.1)

\begin{figure}[here]
\begin{picture}(400,100)(0,0)


\put(100,65){\circle{40}}
\put(300,65){\circle{40}}

\qbezier(115,76)(100,65)(115,54)
\qbezier(85,76)(100,65)(85,54)
\qbezier(90,82)(100,65)(110,82)
\qbezier(90,48)(100,65)(110,48)

\qbezier(115,76)(113,69)(118,70)
\qbezier(118,70)(110,65)(115,54)

\qbezier(85,76)(87,69)(82,70)
\qbezier(82,70)(90,65)(85,54)

\qbezier(90,48)(97,56)(100,45)
\qbezier(100,45)(103,56)(110,48)

\qbezier(90,82)(97,74)(100,85)
\qbezier(100,85)(103,74)(110,82)


\qbezier(320,65)(300,65)(310,80)
\qbezier(310,80)(300,65)(290,80)
\qbezier(290,80)(300,65)(280,65)
\qbezier(280,65)(300,65)(290,50)
\qbezier(290,50)(300,65)(310,50)
\qbezier(310,50)(300,65)(320,65)

\put(30,20){{\sf (i) Maximal number of ideal}}
\put(85,7){{\sf polygons}}
\put(230,20){{\sf (ii) Minimal number of ideal}}
\put(280,7){{\sf polygons}}

\end{picture}
\caption{Region $U=\sqcup ~U_i$ for surface $g=2$}
\end{figure}

\medskip\noindent
The hyperbolic area of an ideal $n$-gon $U_i$ is equal to $(n-2)\pi$,
see \cite{CaB}. Since
\begin{equation}\label{eq4}
\sum Area~U_i=(4g-4)\pi,
\end{equation}
the number of ideal polygons in $U$ is
{\it finite}.

\bigskip
By a  {\it singularity data} we mean  
the number and shape of the ideal polygons $U_i$ which cover the
principal region of $\lambda$. Since
$Area~U_i=(n-2)\pi$ and $\sum~Area~U_i=(4g-4)\pi$, there exists
only finite number of opportunities for such data with fixed $g$.

Let $\Delta=(k_1,\dots,k_m)$ be a set of positive integers
and half-integers such that $\sum k_i=2g-2$.
To each ideal $n$-gon $U_i$ we assign number $k_i$ such that
\begin{equation}\label{eq5}
k_i={n-2\over 2}.
\end{equation}
The reader can verify that condition $\sum~Area~U_i=(4g-4)\pi$
is equivalent to $\sum k_i=2g-2$.
\begin{dfn}\label{dfn1}
Given lamination $\lambda\in\Lambda(S)$, the unordered
tuple $\Delta=(k_1,\dots,k_m)$ is called a singularity
data of $\lambda$.
\end{dfn}

\subsection{Morse coding of the geodesic lines}
The idea of the method is to dissect the surface along 
$r=2g+m-1$ loops so that it  becomes simply connected.
Each loop gets a label and the geodesic line becomes a
bi-infinite sequence of the labels accordingly the order
it intersects the loops. Conversely, every  bi-infinite 
sequence of the labels defines a geodesic (Morse's theorem). 
Let us pass to the construction whose details
can be found in (\cite{Mor}).

Let $S$ be a connected complete  hyperbolic surface of
genus $g$   with $m$ boundary components. With no restriction,
we can assume that the boundary components
\begin{equation}
h_1,h_2,\dots, h_m
\end{equation}
are closed geodesics of $S$. To render $S$ simply connected,
choose a point $P$ on one of the boundary components, say $h_m$. 
First, using geodesic arcs one connects $P$  with a point lying on each of the
remaining boundary components 
\begin{equation}
h_1,h_2,\dots, h_{m-1},
\end{equation}
and then dissects $S$ along the arcs. The new surface will 
have a unique boundary.

Let 
\begin{equation}
c_1,c_2,\dots, c_{2g}
\end{equation}
be the geodesic loops based in the point $P\in h_m$ and such that
they   dissect  $S$ into a simply connected  plane region $T$.

For simplicity, let us assume that the boundary components $p_i$
are punctures. In this case $T$ is a polygon bounded by an
even number of the geodesic segments
\begin{equation}
c_1, c_1', c_2, c_2',\dots, c_{2g}, c_{2g}'; h_1, h_1', h_2, h_2',\dots, h_{m-1}, h_{m-1}',
\end{equation}
where slash denotes the opposite sides of the cut.  

Suppose we have an infinite stock of copies of $T$. If $l$
is a geodesic on $S$ which is disjoint from  the boundary
of $S$, we label each  copy of $T$ with a symbol $\sigma$  from the
set 
\begin{equation}\label{label}
c_1,c_2,\dots, c_{2g}; h_1, h_2,\dots, h_{m-1} 
\end{equation}
if $l\cap\sigma\ne\emptyset$. In this way one constructs 
an infinite cyclic cover of $S$ by gluing the copies of $T$
along the sides of $T$ hit by $l$. Such a cover is uniquely defined 
by a bi-ininite sequence of symbols
\begin{equation}
\dots \sigma_{-1}, ~\sigma_0, ~\sigma_1\dots
\end{equation}
with the values in the set (\ref{label}).  The above bi-infinite
sequence is called a {\it reduced (symbolic) curve}.

Clearly, every geodesic on $S$ which does not intersect the boundary
can be turned into  a ``symbolic curve''.  An amazing fact
proved by M.~Morse is that the converse is true.
\begin{lem}
There is a one to one correspondence between the set of all geodesics
on $S$ which does not intersect the boundary of $S$ and the set of all
reduced curves. 
\end{lem}
{\it Proof.} This is essentially Theorem 3 of (\cite{Mor}).
$\square$

\subsection{$AF$ $C^*$-algebras and Bratteli diagrams}
An $AF$ (approximately finite-dimensional) algebra is defined to
be a  norm closure of an ascending sequence of the finite dimensional
algebras $M_n$'s, where  $M_n$ is an algebra of $n\times n$ matrices
with the entries in $\C$. Here the index $n=(n_1,\dots,n_k)$ represents
a {\it multi-matrix} algebra $M_n=M_{n_1}\oplus\dots\oplus M_{n_k}$.
Let
\begin{equation}\label{e20}
M_1\buildrel\rm\varphi_1\over\longrightarrow M_2
   \buildrel\rm\varphi_2\over\longrightarrow\dots,
\end{equation}
be a chain of algebras and their homomorphisms. A set-theoretic limit
$A=\lim M_n$ has a natural algebraic structure given by the formula
$a_m+b_k\to a+b$; here $a_m\to a,b_k\to b$ for the
sequences $a_m\in M_m,b_k\in M_k$.
The homeomorphisms of the above (multi-matrix) algebras admit a 
canonical description (\cite{E}).
Suppose that
$p,q\in\N$ and $k\in\Z^+$ are such numbers that $kq\le p$. Let us
define a homomorphism $\varphi:M_q\to M_p$ by the formula
\begin{equation}\label{e21}
a\longmapsto\underbrace{a\oplus\dots\oplus a}_k\oplus 0_h,
\end{equation}
where $p=kq+h$. More generally, if $q=(q_1,\dots,q_s),p=(p_1,\dots,p_r)$
are vectors in $\N^s,\N^r$, respectively, and $\Phi=(\phi_{kl})$ is
a $r\times s$ matrix with the entries in $\Z^+$ such that
$\Phi(q)\le p$, then the homomorphism $\varphi$ is defined by
the formula:
\begin{eqnarray}\label{e22}
a_1\oplus\dots\oplus a_s &\longrightarrow&
\underbrace{(a_1\oplus a_1\oplus\dots)}_{\phi_{11}}\oplus
\underbrace{(a_2\oplus a_2\oplus\dots)}_{\phi_{12}}\oplus\dots\oplus 0_{h_1}
\\
&\oplus&  
\underbrace{(a_1\oplus a_1\oplus\dots)}_{\phi_{21}}\oplus
\underbrace{(a_2\oplus a_2\oplus\dots)}_{\phi_{22}}\oplus\dots\oplus 0_{h_2}
\oplus\dots\nonumber
\end{eqnarray}
where $\Phi(q)+h=p$. We say that $\varphi$  is a {\it
canonical homomorphism} between $M_p$ and $M_q$. 
Any homomorphism $\varphi:M_q\to M_p$ can be rendered canonical
(\cite{E}).

Graphical presentation of the canonical homomorphism
is called a {\it Bratteli diagram}. Every ``block''
of such diagram is a bipartite graph with $r\times s$ matrix
$\Phi=(\phi_{kl})$.

In general,
Bratteli diagram is given by a vertex set $V$ and edge set $E$ such that
$V$ is an infinite disjoint union $V_1\sqcup V_2\sqcup\dots$, where each $V_i$ has
cardinality $n$. Any pair $V_{i-1},V_i$ defines a non-empty set
$E_i\subset E$ of  edges with a pair of range and source functions
$r,s$ such that $r(E_i)\subseteq V_i$ and $s(E_i)\subseteq V_{i-1}$.
The non-negative integral matrix of ``incidences''
$M=(\phi_{ij})$  shows how many edges
there are between the $k$-th vertex in row $V_{i-1}$ and $l$-th vertex
in row $V_i$.

\section{Proof of Theorem 1}
Let us outline the main steps of the proof. Let $B$ 
be a strictly ergodic unimodular Bratteli diagram of
rank $r$. Then $B$ defines a simple dimension group $G$ 
of rank $r$ with a unique state. Any such group can
be realized as a dense subgroup ${\Bbb Z}\lambda_1+\dots
+{\Bbb Z}\lambda_r$ of the real line. For a canonical state, we have
$\sum \lambda_i=1$. Given the singularity 
data $\Delta$, we can construct an interval exchange
transformation $(\lambda,\pi)$ on the intervals $\lambda_1,
\dots,\lambda_r$. The infinite sequence of induced
transformations $(\lambda,\pi)\supset (\lambda',\pi')
\supset\dots$ contracts to a point $\theta\in [0,1]$.    
It is possible to associate to $\theta$ an infinite
sequence of symbols taking values in the set
$\lambda_1,\dots,\lambda_r$. We show that such a sequence
is a recurrent non-perodic symbolic geodesic sequence. 
By the Morse theorem, we get a recurrent geodesic on
the surface of genus $g$. Let us pass to the detailed
construction.

\bigskip
{\sf Part I.}
Let $(B,\Delta)$ be as in the theorem. We wish to construct 
an interval exchange transformation $(\lambda,\pi)$ from
the pair $(B,\Delta)$. For that consider a dimension group 
\begin{equation}
G(B)=\lim_{n\to\infty}({\Bbb Z}^r, M_n),
\end{equation}
where $r$ is the rank of $B$ and $M_n$ are matrices of partial
multiplicities of $B$. The unimodularity of $B$ implies that $G$
is a simple dimension group of the rank $r$ without infinitesimal elements. 
The strict ergodicity of $B$ is equivalent to the state space $S(G)$ of $G$  is a point (\cite{E},
Ch. 4). Let us recall the following lemma.
\begin{lem}\label{lm3}
Suppose that $G$ is a simple dimension group of rank $r$
without infinitesimal elements and $dim~S(G)=d-1$. Then 
$G$ is order isomorphic to a dense subgroup of ${\Bbb R}^d$, provided
with the relative strict order.  
\end{lem}
{\it Proof.}  This is essentially a Corollary 4.7 p. 25 of  (\cite{E}). 
$\square$

\medskip\noindent
Note that in our case $dim~S(G)=0$ and therefore $d=1$. Thus
by Lemma \ref{lm3} we have a dense subgroup of rank $r$ of
the real line. Let us fix generators  $\lambda_1,\dots,\lambda_r$
of the subgroup to be positive reals. Note that $\lambda_i$ are
linearly independent over ${\Bbb Q}$ except the normalization
condition $\lambda_1+\dots+\lambda_r=1$ which comes from the 
unique standard state on $G$.

We set $\lambda=(\lambda_1,\dots,\lambda_r)$ and we wish to 
construct  a permutation $\pi$ on the above intervals from the 
singularity data $\Delta$. Recall  that every element $\pi\in \Sigma_r$
of the permutation group on $r$ elements decomposes into the
elementary cycles $\pi_1\circ\dots\circ\pi_s$. The decomposition
is unique up to a cyclic permutation.

Let $\Delta=(k_1,\dots,k_m)$ be a singularity data. The Veech's 
``zippered rectangles'' construction (Section 6 of \cite{Vee2})
implies that the total number of elementary cycles 
\begin{equation}
s=r-2g+1=m.
\end{equation}
The length of the elementary cycle $\pi_i$ is also determined
(up to an isomorphism) by the corresponding singularity $k_i$,
see Veech, {\it ibid.} Therefore, we get a permutation $\pi\in\Sigma_r$
such that
\begin{equation}
\pi= \pi_1\circ\dots\circ\pi_m.
\end{equation}

\bigskip
{\sf Part II.}
Let $(\lambda,\pi)$ be an interval exchange transformation 
obtained from the pair $(B,\Delta)$. We will assign to $(\lambda,\pi)$
an infinite sequence of symbols taking value in the finite set $\lambda_1,
\dots,\lambda_r$. To achieve this goal, we  will use  the concept of ``induced    
transformations'' developed  by Keane, Rauzy and Veech (\cite{Vee1}).

Recall that an interval 
\begin{equation}
\Gamma=[\xi,\eta), \qquad 0\le\xi<\eta\le |\lambda|,
\end{equation}
where $\lambda=\lambda_1+\dots+\lambda_r$, is called {\it admissible}
for the interval exchange transformation $\varphi=\varphi(\lambda,\pi)$ 
if $\Gamma$ splits on $r$ parts $\lambda'_1,\dots,\lambda'_r$ 
such that $\varphi^{n(\lambda'_i)}$ is continuous on each of $\lambda'_i$. 
The corresponding interval exchange transformation on $\Gamma$ is
called {\it induced}. 
The positive vectors $\lambda=(\lambda_1,\dots,\lambda_r)$ and 
$\lambda'=(\lambda'_1,\dots,\lambda'_r)$ are connected by the
formula
\begin{equation}
\lambda=M\lambda',
\end{equation}
where $M$ is a non-negative integral matrix of determinant $\pm 1$,
see (\cite{Vee1}, Section 3). The matrix $M$ coincides with the matrix
of the partial multiplicity $M_1$ which occurs at the first position in the Bratteli
diagram  $B$.

Let 
\begin{equation}
\Gamma_1\supset\Gamma_2\supset\dots 
\end{equation}
be an infinite sequence of admissible intervals. Clearly, $|\Gamma_n|\to 0$
as $n\to\infty$. Then the set 
\begin{equation}
\theta=\bigcap_{n=1}^{\infty}\Gamma_n,
\end{equation}
is either empty or conists of a point. Assuming that the admissible intervals
have only a finite number of  common (right) endpoints, we get that $\theta$
is a point such that $0<\theta<1$.

Denote by $\lambda_i^{(j)}, ~1\le i\le r$ a part of admissible interval $\Gamma_j$
such that $\theta\in \lambda_i^{(j)}$.  A sequence of subintervals
\begin{equation}
S=\left\{\lambda_{i_j}^{(j)}\right\}_{j=1}^{\infty},
\end{equation}
we call a {\it pre-code}.

To construct a {\it code} $S^*$ from the pre-code $S$, we insert a finite
number of symbols between any two symbols $\lambda_{i_{j-1}}^{(j-1)}$
and $\lambda_{i_j}^{(j)}$ of $S$ as follows.      
Recall that
\begin{equation}
\lambda_{i_{j-1}}^{(j-1)}= a_{i_j 1}\lambda_1^j+\dots+a_{i_j r}\lambda_r^j,
\end{equation}
where $a_{mn}$ are entries of the matrix $M_j$.  We insert  $a_{i_j 1}$ symbols
$\lambda_1^j$,  $a_{i_j 2}$ symbols $\lambda_2^j$, etc,  between the symbols 
 $\lambda_{i_{j-1}}^{(j-1)}$ and $\lambda_{i_j}^{(j)}$ of $S$ in the order the
orbit of the point $\theta$ under the induced transformation $\varphi_j=\varphi_j(\lambda_j,\pi_j)$
hits the admissible  interval $\Gamma_j$. We have therefore:
\begin{equation}
\dots \lambda_{i_{j-1}}^{(j-1)}\quad \underbrace{\lambda_1^j\dots\lambda_1^j}_{a_{i_j 1}}
\dots
\underbrace{\lambda_r^j\dots\lambda_r^j}_{a_{i_j r}}
\quad\lambda_{i_j}^{(j)} \dots
\end{equation}

\bigskip
{\sf Part III.}
Let $S^*$ be a code associated to the pair $(B,\Delta)$ as described above. 
$S^*$ can be converted to a {\it symbolic geodesic} $\Sigma=\Sigma (S^*)$ 
by ``forgetting'' the upper indices  in the sequence $S^*$. Clearly, the symbols
of $\Sigma$ take values in a finite set of cardinality  $r$. 
\begin{lem}
$\Sigma$ is a recurrent non-periodic symbolic geodesic. 
\end{lem}
{\it Proof.} The idea is to identify $\Sigma$ with a recurrent
trajectory of a suspension flow over the interval exchange 
transformation $\varphi$ constructed in Part I. Indeed,
let $v_t$ be such a flow obtained by  the ``zippered rectangles'' method (\cite{Vee2}).
Consider a trajectory $l=v_t(\theta)$ through the point $\theta$ 
defined in Part II. Since flow $v_t$ is minimal, the closure of $l$
is the entire surface $S$. In particular, $l$ is a recurrent 
non-periodic trajectory. 

The intervals 
\begin{equation}
\lambda_1,\dots,\lambda_r
\end{equation}
give a dissection of $S$ into a simply connected
domain as follows.  For $i=1,\dots,r$,   one takes a rectangle with the
opposite sides $\lambda_i, \lambda_i'$  and $f_t(a),
f_t(b)$, where $\lambda_i'$ is the image of $\lambda_i$
under the Poincar\'e (first return) mapping
and $a,b$ are the ends of the interval $\lambda_i$.     
In this way, the recurrent trajectory $l$ becomes
a symbolic trajectory $\Sigma$ with the desired property. 
 $\square$

\bigskip
To finish the proof, let $S$ be a hyperbolic surface.
Take a standard dissection of $S$ by the $r$ curves 
$c_1,\dots,c_{2g}; h_1,\dots, h_{m-1}$ as described in Section
1.2. Then the Morse theorem says that there exists a recurrent
non-periodic geodesic $l_{\Sigma}$ on $S$. The closure
$\lambda=\overline{l_{\Sigma}}$ is a minimal geodesic 
lamination with $m$ principal regions.  Theorem follows. 
$\square$

\section{An example}
In this section we consider an example of the ``golden mean''
Bratteli diagram. We construct a symbolic geodesic $\Sigma$
in this case, and show that $\Sigma$ coincides with an example
of Morse. 
 
\bigskip\noindent
{\sf Example.}
Let $B$ be a Bratteli diagram presented in Fig. 2.

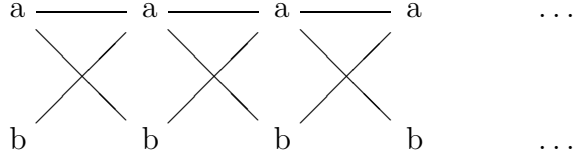
\begin{figure}[here]
\begin{picture}(400,100)(0,0)

\put(50,75){a}
\put(50,25){b}

\put(100,75){a}
\put(100,25){b}

\put(150,75){a}
\put(150,25){b}

\put(200,75){a}
\put(200,25){b}

\put(60,77){\line(1,0){34}}
\put(110,77){\line(1,0){34}}
\put(160,77){\line(1,0){34}}

\put(60,70){\line(1,-1){34}}
\put(110,70){\line(1,-1){34}}
\put(160,70){\line(1,-1){34}}

\put(60,35){\line(1,1){34}}
\put(110,35){\line(1,1){34}}
\put(160,35){\line(1,1){34}}

\put(250,75){$\dots$}
\put(250,25){$\dots$}

\end{picture}
\caption{Golden mean diagram.}
\end{figure}

\medskip\noindent
The incidence matrix is a constant unimodular matrix
\begin{equation}
M_n=M=\left(\matrix{1 &1\cr 1 & 0}\right). 
\end{equation}
Such a diagram is known to be strictly ergodic, see Effros (\cite{E}), Theorem 6.1.   
The singularity data $\Delta=(0)$, i.e. there is a unique singular point 
of index $O$ (a fake saddle). 

We wish to construct a recurrent geodesic $l_{\Sigma}$ from the pair $(B,\Delta)$. 
First, notice that the formula
\begin{equation}
r=2g+m-1,
\end{equation}
 implies $g=1$ since $r=2, m=1$. Therefore, our surface is a torus $T^2$. 

The canonical state on the dimension group $G=G(B)$ gives us 
\begin{equation}
\lambda_1= {\sqrt{5}-1\over 2},\qquad \lambda_2={3-\sqrt{5}\over 2}.
\end{equation}

The sequence of admissible intervals $\Gamma_1\supset\Gamma_2\supset\dots$ becomes 
\begin{equation}
[\lambda_1-{\lambda_1\over\varepsilon},\quad\lambda_1+{\lambda_2\over\varepsilon}]\supset
[\lambda_1-{\lambda_1\over\varepsilon^2},\quad\lambda_1+{\lambda_2\over\varepsilon^2}]\supset\dots
 \end{equation}
where $\varepsilon={3+\sqrt{5}\over 2}$ is the Perron-Frobenius eigenvalue of the matrix
$M^2$. Such a sequence contracts to the point $\theta=\lambda_1={\sqrt{5}-1\over 2}$.

Denote by $a$ and $b$ the points of the upper and lower row of the Bratteli diagram on Fig.2. 
Then the pre-code of $\theta$:
\begin{equation}
S=b\underbrace{a}_1 b\underbrace{aa}_2 b\underbrace{aaa}_3b\dots 
\end{equation}

The code of $\theta$ is obtained from $S$ by inserting $a$ between any $ba$, $b$ between any $aa$ and
nothing between any $ab$: 
\begin{equation}
S^*=b{\bf a}\underbrace{a}_1 b{\bf a}\underbrace{a{\bf b}a}_2 b\underbrace{a{\bf b}a{\bf b}a}_3b\dots
\end{equation}

Up to a  notation, $S^*$ coincides with the Morse example of a forward symbolic  non-periodic recurrent
geodesic on $T^2$, see (\cite{Mor}), \S 14.


\end{document}